
\def\calL{{\cal L}}
\def\calH{{\cal H}}

\def\calK{{\cal K}}

\def\calO{{\cal O}}

\def\char{{\rm char\ }}
\def\deg{{\rm deg}}

\font\tenfrak=eufm10
\font\sevenfrak=eufm7
\font\fivefrak=eufm5
\newfam\frakfam
\textfont\frakfam=\tenfrak
\scriptfont\frakfam=\sevenfrak
\scriptscriptfont\frakfam=\fivefrak
\def\frak{\fam\frakfam\tenfrak}

\font\tenbb=msbm10
\font\sevenbb=msbm7
\font\fivebb=msbm5
\newfam\bbfam
\textfont\bbfam=\tenbb
\scriptfont\bbfam=\sevenbb
\scriptscriptfont\bbfam=\fivebb
\def\bb{\fam\bbfam\tenbb}
\font\tensym=msam10
\font\sevensym=msam7
\font\fivesym=msam5
\newfam\symfam
\textfont\symfam=\tensym
\scriptfont\symfam=\sevensym
\scriptscriptfont\symfam=\fivesym

\font\tenlas=lasy10
\font\sevenlas=lasy7
\font\fivelas=lasy5
\newfam\lasfam
\textfont\lasfam=\tenlas
\scriptfont\lasfam=\sevenlas
\scriptscriptfont\lasfam=\fivelas

\def\bbZ{{\bb Z}}

\def\bbC{{\bb C}}
\def\bbQ{{\bb Q}}
\def\bbP{{\bb P}}

\def\bbA{{\bb A}}

\def\calF{{\cal F}}

\def\calH{{\cal H}}
\def\calL{{\cal L}}

\def\calN{{\cal N}}
\def\calO{{\cal O}}

\def\calK{{\cal K}}
\def\calE{{\cal E}}
\def\calA{{\cal A}}
\def\Ker{{\rm Ker }}
\def\Pic{{\rm Pic}}
\def\Ima{{\rm Im }}
\def\Spec{{\rm Spec}}
\def\Hom{{\rm Hom}}
\def\small{\medskip\noindent}
\def\med{\medskip\noindent}

\magnification = 1100
\vglue .3 in
\centerline {{\bf WILD p-CYCLIC ACTIONS  ON
K3-SURFACES}}
\vglue .4in
\centerline{ I. DOLGACHEV \footnote{*}{
Research supported in part by a NSF grant.} and J.
KEUM\footnote{$^{**}$} {Supported by Korean Ministry of  Education(BSRI
1438) and KOSEF-GARC.}}

\vglue .3 in
\med
\centerline {{\bf  0. Introduction}}
\med
\med
 In this paper we study automorphisms $g$ of order
$p$ of K3-surfaces defined over an algebraically closed field of characteristic $p > 0$. We divide all
possible actions in the following cases according to the structure of the 
set of fixed points $X^g$:
$X^g$ is a finite set, $X^g$ contains a one-dimensional part $D$ which is a positive divisor of Kodaira
dimension
$\kappa(X,D) = 0,1,2$. In the latter case we prove that $X^g = D$ and $D$
is connected. Accordingly we obtain the following results. 

\med
\proclaim Theorem 1.  Suppose that the fixed locus $X^g$ of $g$
is finite. Then $|X^g|\le 2$. If $X^g = \emptyset$, then $p = 2$ and $Y = X/(g)$ is a
non-classical  Enriques surface of $\mu_2$-type and $X$ is its K3-cover. If 
$|X^g| = 1$ then a minimal resolution $\tilde Y$ of $Y$ is either a rational surface
or a K3-surface. If
$|X^g| = 2$, then
$\tilde Y$ is a K3-surface. The cases when $\tilde Y$ is a
K3-surface can occur only if
$p
\le 5$.

\med
Suppose now that $X^g$ is one-dimensional. The next results cover
three different cases corresponding to possible values of the Kodaira
dimension (or $D$-dimension) $\kappa(X,X^g)$ of the divisor $X^g$ of fixed
points.

\med
\proclaim Theorem 2.  Suppose that 
$\kappa(X,X^g) = 0$. 
 Then $X^g$ is a nodal cycle, i.e. the union of smooth
rational curves whose intersection matrix is negative definite 
(of type A-D-E). Moreover, if we contract
$X^g$ to get a singular normal surface $X'$, then  the quotient $X'/(g)$ is a
rational surface with one rational singular point or one  elliptic singular
point. If $p\ge 5$ then the former case does not occur.

\med
\proclaim Theorem 3.  Suppose that $\kappa(X,X^g) = 1$. 
 Then $p\le 11$ and we have the following: 
\item{(1)} There exists a divisor $D$ with support equal to $X^g$ such that
the linear system $|D|$  defines an elliptic or quasi-elliptic fibration 
$\phi : X \to \bbP^1$ (the second option can occur only if $p \le 3$).
\item{(2)} The induced map $g^*$ on the base curve $\bbP^1$ is of
order
$p$ unless $p = 2$ and $\phi$ is an elliptic fibration.
\item{(3)} If $p > 3$ and $\phi$ is an elliptic fibration, the type of the fibre $D$
in Kodaira's notation is 
$II$ if $p=11$, $III$ or $II^*$ if $p=7$, and 
$IV$ or
$III^*$ if
$p=5$. 
\item{(4)} If $\phi$ is a quasi-elliptic fibration, the type of the fibre $D$ is 
$IV$ or $II^*$ if $p = 3$, and of type $II, III, II^*, I_{2k}^*$, $0\le k \le 8$ if $p
= 2$.

\med
Note that the assumption $p\ne 2$ in (2) is esential as we show by an example.

\med
\proclaim Theorem 4. Suppose that $\kappa(X,X^g) = 2$. 
 Then $X^g$ is equal to the support of some nef divisor
$D$ with $D^2 > 0$. Assume that it is chosen to be minimal with this
property. Let $d = \dim H^0(X,\calO_X(D-X^g))$ and $N = {1\over 2}D^2+1$.
Then
$p(N-d-1)\le 2N-2$.
In particular, if $X^g$ is irreducible and $p\ne 2$, then the pair $(p,N) =
(3,2),(3,3),(3,4)$, or $(s,2), s\ge 5$. Moreover, if $p\ne 2$, $X^g$
is a singular irreducble  rational curve. 

\med
The last theorem is essentially a result from the earlier work of 
Oguiso [Og2] where it was assumed that $X^g$ is irreducible. In this case  explicit
computations show that $s = 5$. Oguiso also
obtains the defining equation of
$X$ in a suitable weighted projective space. We note here that 
Oguiso's method does not
apply to other cases, as it uses the projective model of $X$ determined by the linear
system $|X^g|$ on which $g$ acts linearly.

Note that if an automorphism $g$ of a K3-surface over $k$ with
positive characteristic
 $p$
has an order of a power of $p$, then $g$ is symplectic, i.e. acts on the space of global
 $2$-forms trivially. Over the field of complex numbers $\bbC$, Nikulin
([Ni1], \S 5)  determined
all finite abelian symplectic automorphism groups of K3-surfaces. In particular, he
proved that if $g$ is a symplectic automorphism of order $n\ge2$ of a K3-surface over 
\bbC, then the fixed locus $X^g$ is nonempty and finite and the pair $(n,|X^g|)$ is
one of the following : $(2,8),(3,6),(4,4),(5,4),(6,2),(7,3)$ and $(8,2)$. This result, 
being deduced from the fact that the quotient surface $X/(g)$ becomes a K3-surface 
with only
rational double points of type $A_n$, can be extended to the case of positive characteristic
$p$ as long as $p\ge 3$ and the order $n$ of the automorphism is coprime to $p$.

Also for comparison, note that involutions $g$ of prime order $p$ on a
K3-surface over a field of characteristic prime to $p$ with
non-isolated fixed points can be completely classified. We have $p = 2$ or
$3$. If $p = 2$ the classification can be found in [Ni2]. If $p =3$,
then $X/(g)$ is a nonsingular quadric surface and the map $\pi:X\to X/(g)$ is
branched along a smooth curve of bidegree $(3,3)$. This easily follows
from the Hurwitz formula for the cyclic cover $\pi$ of nonsingular
surfaces. 

In section 1 we give some general results about wild cyclic actions, and, in particular, answer a
question from [Pe1]. In section 2 we study the case of isolated fixed
points where our main tool is Grothendieck's spectral sequence for
cohomology of $G$-equivariant sheaves from [Gr]. Section 3 is devoted to
the case of non-isolated fixed points where we establish some general
results in this case, proving for example that the fixed locus must be
connected. We also show that the quotient surface $Y = X/(g)$ is a normal
anticanonical surface with rational singularities. The other sections 4-6
are devoted to more detailed study of three different cases corresponding
to the possible Kodaira dimension of the divisor of fixed points. At the end of
section 5 we give some applications of our results to the theory of elliptic
fibrations on a K3-surface.

A part of this paper was written during the second author's stay at
Nagoya University June of 1998 on JSPS fellowship. He thanks S. Kond$\bar {\rm o}$
and the University for hospitality, and JSPS for financial support.
He also thanks K. Oguiso for bringing to his attention the paper [Og1] which contains the
key idea for the proof of Lemma 5.2. Discussions with S. Mukai and S. Mori
are very much appreciated.
\bigskip

\centerline{{\bf 1.  Preliminaries on wild cyclic actions}}.
\med

Let $G$ be a group acting on a topological space $X$, and let $Y=X/G$ be
the  quotient space and $\pi : X \to Y$ the quotient map. Consider the category ${\cal
S}(X,G)$ of abelian $G$-sheaves on $X$ and the 
category  ${\cal A}$ of abelian groups. The functor
$${\cal S}(X,G) \to {\cal A},\quad {\cal F} \mapsto \Gamma(X, {\cal F})^G$$
can be represented as a composition of functors in two different ways:
$${\cal S}(X,G) \to {\cal S}(Y) \to {\cal A},\quad {\cal F} \mapsto \pi_*^G{\cal F}
\mapsto
 \Gamma(Y, \pi_*^G{\cal F}),$$
$${\cal S}(X,G) \to {\cal A} \to {\cal A},\quad {\cal F}  \mapsto \Gamma(X, {\cal F})
\mapsto \Gamma(X, {\cal F})^G,$$
where  $\pi_*^G{\cal F}$ is the sheaf $U \mapsto  \Gamma(\pi^{-1}(U),{\cal F})^G$. This 
gives two spectral sequences for the compositions of functors (see [Gr],
Theorem 5.2.1):
$$'E_2^{p,q}(\calF) = H^{p}(Y,{\cal H}^{q}(G,{\cal F})) 
\Rightarrow H^n,\eqno (1.1)$$
$$''E_2^{p,q}(\calF) = H^{p}(G,H^{q}(X,{\cal F})) \Rightarrow H^n.\eqno (1.2)$$
Here the sheaves ${\cal H}^{q}(G,{\cal F})$ can be defined on open
subsets of $Y$ by
$${\cal H}^{q}(G,{\cal F})(U) =
H^q(G,\calF|\pi^{-1}(U)).\eqno (1.3)$$
 We will
apply this to the situation when
$X$ is an irreducible algebraic variety over a field $k$, 
$G$ is a finite group of its automorphisms and
${\cal F}$ is the structure sheaf ${\cal O}_X$.
\med

\proclaim Proposition 1.1.  Let $G$ be a finite group acting freely on an
irreducible algebraic variety $X$ over a field $k$.  Then
$$\calH^{i}(G,\calO_X) = 0, i > 0.$$

{\sl Proof.} Let $\pi:X\to Y$ be the quotient map. It is a finite
and hence affine morphism. So, without loss of generality we may assume
that
$X = \Spec(A)$
is an $k$-algebra and $Y = \Spec(B)$, where $B = A^G$ is the subalgebra of
invariant elements. We need to show that
$H^i(G,A) = 0$ for
$i > 0$. Since $G$ acts freely, $X$ is a principal fibre bundle over
$Y$ in the sense that $A$ is flat over $A^G$ and the canonical
morphism $G\times X\to X\times_YX$ defined by the action is an
isomorphism (see [Mu], Proposition 0.9). Consider the group
algebra $B[G]$ of $G$ over $B$ as a $B$-module. By the above we have an
isomorphism of $A$-modules 
$$B[G]\otimes_BA \cong (k[G]\otimes_kB)\otimes_BA \cong k[G]\otimes_kA
\cong A\otimes_BA.$$
Since $A$ is flat over $B$, this implies that $A$ is isomorphic to
$B[G]$ as a $B$-module. As a $G$-module, $A$ is induced from the trivial
$\{1\}$-module $B $. By Shapiro's Lemma, we get
$$H^i(G,A) \cong H^i(\{1\},B) = 0,\quad i > 0.$$
This proves the assertion. 

\med
\proclaim Corollary 1.2. Let $G$ be a finite group acting on an
irreducible variety $X$. Then the sheaves $\calH^q(G,\calO_X), q > 0,$
are equal to zero over the quotient of the open subset of $X$ where $G$
acts freely.

\med
Now let $G = (g)$ be a cyclic group of order $n$.  Recall that the cohomology of a cyclic 
group $G = (g)$ of order $n$ with coefficient in a $G$-module 
$M$ are computed as follows:
$$H^{2i+1}(G,M) = \Ker T/\Ima(g-1),\quad H^{2i}(G,M) = \Ker(g-1)/\Ima T, i > 0, \quad
H^0(G,M) = \Ker(g-1),$$ where 
$T  = 1+g+\ldots+g^{n-1}$ and $g-1$ are elements of the group algebra $\bbZ[G]$ acting 
naturally on $M$.
Globalizing we get the homomorphisms of $\calO_Y$-modules
$$T:\pi_*(\calO_X)\to \calO_Y,\quad g-1:\pi_*(\calO_X)\to \pi_*(\calO_X)$$
and the isomorphisms of $\calO_Y$-modules
$$\calH^{2i+1}(G,\calO_X) = \Ker T/\Ima(g-1), \quad \calH^{2i}(G,\calO_X) =
\Ker(g-1)/\Ima T,\quad  i > 0.\eqno (1.4)$$
Of course, $\calH^0(G,\calO_X) = \calO_Y$.

\med
\proclaim Proposition 1.3. Let $G = (g)$ be a cyclic group acting on a
Cohen-Macaulay algebraic variety
$X$. Assume that the quotient $Y = X/G$ is also Cohen-Macaulay (e.g. $\dim X = 2$). Then we 
have an exact sequence of sheaves on $Y$:
$$0\to \omega_Y \to (\pi_*\omega_X)^G\to 
{\cal E}xt^1_{\calO_Y}(\calH^2(G,\calO_X),\omega_Y)\to 0$$
where $\omega_X$ and $\omega_Y$ are the canonical dualizing sheaves of $X$ and $Y$, 
respectively.

{\sl Proof.} Let $\phi:\omega_Y\to (\pi_*\omega_X)^G$ be the natural homomorphism induced by 
the inverse image of a differential. This map is injective since the map $\pi:X\to Y$ is 
separable.
By  Grothendieck's duality theory (see [Ha], V.2.4),
$$\pi_*(\omega_X) = {\cal H}om_{\calO_Y}(\pi_*(\calO_X),\omega_Y)\eqno
(1.5)$$ and the image of $\phi$ is the subsheaf of $G$-invariant
$\calO_Y$-module  homomorphisms which are composed of the trace map
$T:\pi_*(\calO_X)\to
\calO_Y$ and a map $\calO_Y\to
\omega_Y$ (see [Pe3]). Consider the natural exact sequence
$$0\to \Ker(T)/\Ima(g-1)\to \pi_*(\calO_X)/\Ima(g-1)\to  \Ima(T)\to 0. \eqno (1.6)$$
Clearly, the dual sheaf of $\pi_*(\calO_X)/\Ima(g-1)$ is equal to the 
subsheaf of $G$-invariant $\calO_Y$-module 
homomorphisms $\pi_*(\calO_X)\to \omega_Y$. Since the sheaf
$\calH^1(G,\calO_X) =\Ker(T)/\Ima(g-1)$ is a torsion $\calO_Y$-module,  we
get an isomorphism
$${\cal H}om_{\calO_Y}(\Ima(T),\omega_Y)\cong
{\cal H}om_{\calO_Y}(\pi_*(\calO_X)/\Ima(g-1),\omega_Y).$$
Using the obvious exact sequence
$$0\to {\cal H}om_{\calO_Y}(\calO_Y,\omega_Y) \to {\cal
H}om_{\calO_Y}(\Ima(T),\omega_Y)
\to {\cal
E}xt^1_{\calO_Y}(\calO_Y/\Ima(T),\omega_Y)  \to 0$$  and the fact that the map $\phi$
is the dual of the trace homomorphism $T$ ([Pe3], Proposition (2.4)), we get the
exact sequence from the assertion of the proposition. 

\med
\proclaim Corollary 1.4 (see [Pe3], Thm.2.7). Assume additionally that $g$ acts freely 
outside a closed subset of codimension $\ge 2$ or
$n$ is invertible in $k$. Then
$$\omega_Y \cong \pi_*(\omega_X)^G.$$

{\sl Proof.} If $n$ is invertible in $k$, the trace homomorphism is surjective so,
$\calH^2(G,\calO_X) = 0$. If $g$ acts freely outside a closed subset of codimension
$\ge 2$ then  $\calH^2(G,\calO_X) \cong\calO_Z$ for some closed subscheme
$Z$ of codimension $\ge 2$.  Therefore, by the property of the dualizing sheaf, we have
${\cal E}xt^1(\calO_Z,\omega_Y) = 0$.

\med
The next proposition answers a question from [Pe1]:
\med
\proclaim Proposition 1.5.  
$$\chi(Y,\calH^1(G,\calO_X)) = \chi(Y,\calH^2(G,\calO_X)).$$

{\sl Proof.} We shall consider the operators $T$ and $g-1$ as endomorphisms of the
$\calO_Y$-Module
$\pi_*(\calO_X)$. Set
$$\calK_i = \Ker(g-1)^{p-i}/\Ker(g-1)^{p-i}\cap \Ima(g-1),\quad \calL_i =
\calO_Y/\Ima(g-1)^{p-i}
\cap \calO_Y.$$ Obviously,
$$\calK_1 = \calH^1(G,\calO_X),\quad \calL_1 = \calH^2(G,\calO_X),$$
$$\calK_{p-1} =\calL_{p-1} = \calO_Y/\Ima(g-1)\cap \calO_Y.$$
We claim the existence of the exact sequence
$$0\to \calK_{i+1}\to \calK_i\to \calL_i\to \calL_{i+1}\to 0,\eqno(1.7)$$
where the map
$\phi_i:\calK_{i}\to \calL_i$
is defined by applying the operator $(g-1)^{p-i-1}$ and the other maps are the natural 
inclusion
$\calK_{i+1}\to \calK_i$ and the surjection $\calL_i\to \calL_{i+1}$. To check the exactedness,
we may assume that $X = \Spec(A), Y = \Spec(A^G)$ are affine. Let $K_i = \calK_i(Y), L_i
=\calL_i(Y).$ If
$a\in \Ker(g-1)^{p-i}$ is a representative of an element from $K_i$, then 
$(g-1)((g-1)^{p-i-1}a) =(g-1)^{p-i}(a) = 0$. Therefore $(g-1)^{p-i-1}(a)\in A^G = B$. Also, 
if $a\in (g-1)(A)\cap \Ker(g-1)^{p-i}$, then $(g-1)^{p-i-1}(a) \in (g-1)^{p-i}(A)$, so
the map is well-defined. 

Let us find its kernel. In the above notation, suppose
$(g-1)^{p-i-1}(a) = (g-1)^{p-i}(x)$ for some $x\in A$. Then
$(g-1)^{p-i-1}(a-(g-1)x) = 0.$ Replacing $a$ with $a-(g-1)(x)$ we may assume that
$a\in \Ker(g-1)^{p-i-1}$ and hence defines an element of $K_{i+1}$. This gives us
a  homomorphism $\Ker(\phi_i)\to K_{i+1}.$ It is clearly bijective.  Since
$(g-1)^{p-i-1}(\Ker(g-1)^{p-i}) = (g-1)^{p-i-1}(A)\cap B$,
the cokernel of $\phi_i$ is equal to $L_{i+1}$.

Now the assertion follows easily from the exact sequence (1.7). We have 
$$\chi(Y,\calK_{i}) -
\chi(Y,\calL_{i}) = \chi(Y,\calK_{i+1}) - \chi(Y,\calL_{i+1}).$$
Starting from $i = p-1$ and ``going down'', we get
$$\chi(Y,\calH^1(G,\calO_X)) - \chi(Y,\calH^2(G,\calO_X)) =
\chi(Y,\calK_{1}) -
\chi(Y,\calL_{1}) = \chi(Y,\calK_{2}) - \chi(Y,\calL_{2}) = 0.$$

\med
Now let us remind some facts about the Artin-Schreier coverings of
algebraic surfaces (see [Ta]).

Recall first a well-known fact from algebra: any cyclic extension of
degree $p$ of a field $K$ of characteristic $p > 0$ is the decomposition
field of an equation $x^p-x = a$ for some $a\in K$. We are trying to
globalize this fact. 

Let $\pi:X\to Y = X/G$ be as in Proposition 1.3 with $G = (g)$ of order $p =\char k$.

\proclaim Proposition 1.6. There is a canonical filtration of 
$\calO_Y$-modules for $\pi_*(\calO_X)$:
$$\pi_*(\calO_X) = \calF_{p-1} \supset \ldots \calF_{1}\supset
\calF_{0} =
\calO_Y,$$
whose quotients $\calL_i = \calF_i/\calF_{i-1}, i = 1,\ldots,p-1,$
are ideal sheaves in $\calO_Y$ satisfying the inclusions:
$$\calL_i\subset \calL_1, \quad \calL_1^i\subset \calL_i.$$ 
If $Y$ is nonsingular, all
$\calF_i$ and
$\calL_1$ are locally free.

{\sl Proof.} Allmost all the assertions are proven in [Ta]. For completeness sake
we briefly remind his
proof. First we may assume that $Y$ and $X$ are affine, say $Y = \Spec(B), X = \Spec(A)$. The 
group algebra $k[G]\cong
k[t]/(t^p-t)$ and the coaction of $G$ on $A$ is given by a homomorphism $\sigma:A\to k[G]
\otimes A = A[t]/(t^p-t).$ For any
$a\in A$ we can write
$$\sigma(a) = a_0+a_1t+\ldots+a_{p-1}t^{p-1}$$
and the axioms of the action imply that the homomorphisms 
$\delta_i:A\to A, \quad a\to a_i,$
satisfy $i!\delta_i = \delta_1^i, \delta_0 = {\rm id}_A$. By definition of the action
$$g(a) = a+\delta_1(a)+\ldots+\delta_{p-1}(a).$$
We set
$$F_i = \{a\in A: \delta_{i+1}(a) = 0\},\quad L_i = \delta_i(F_i).$$
Since, for any $a\in F_i$,  
$$\sigma(\delta_i(a)) = \delta_i(a)+\delta_1(\delta_i(a))+\ldots+\delta_{p-1}(\delta_i(a)) = 
\delta_i(a),$$
we see that $L_i$ is an ideal in $B = A^G$. Moreover, we have an isomorphism of $B$-modules
$$\delta_{i}: F_i/F_{i-1}\to L_i.$$
After globalizing we get the filtration $(\calF_i)$ and the corresponding 
quotient ideal sheaves $\calL_i$. We refer to
[Ta] for the proof of locally freeness of $\calF_i$ and $\calL_1$ in the case when $Y$ is 
nonsingular. It
remains to check the inclusions for the ideal sheaves. Again it is enough to consider the 
affine case.
The first inclusion 
$L_i \subset L_1$ is obvious. By induction, it suffices to show that $L_1L_{i-1}\subset L_i$. 
Let $x\in F_1, y\in F_{i-1}$,
we have
$$g(xy) = g(x)g(y) = (x+\delta_1(x))(y+\delta_1(y)+\ldots+\delta_{i-1}(y)) = \delta_1(x)
\delta_{i-1}(y)+z\in F_{i}$$
where $z\in F_{i-1}$. Clearly,
$\delta_1(x)\delta_{i-1}(y) = \delta_{i}(xy) \in L_{i}$ proving the assertion.
\med

\proclaim Corollary 1.7. Outside a finite set $S$ of points in $Y$ the sheaves $\calL_i$ 
are locally free and $$\calL_i = \calL_1^{\otimes i}.$$
There exists an open affine covering $\{U_i\}$ of $Y\setminus S$ such that
$$\pi^{-1}(U_i) \cong \Spec(\calO_Y(U_i)[t_i]/(t_i^p-a_it_i-b_i))$$
for some $a_i,b_i\in \calO_Y(U_i).$ The elements $a_i$ define a section of $\calL_1^{-p+1}$ 
which is a $(p-1)$-th power of a
section $s$ of $\calL_1^{-1}$. The group $G$ acts on $X$ by the formula $t_i\to t_i+s_i$, 
where $s_i = s|U_i$.

{\sl Proof.} The first assertion is obvious, just take the open subset of the smooth locus of 
$Y$ where the primary components of the ideal sheaves
$\calL_i$ are all of codimension 1. For the remaining statements see [Ta], Lemma 1.2. 

\med
We set 
$$\calL = (\calL_1)^{**} = {\cal H}om_{\calO_Y}({\cal H}om_{\calO_Y}(\calL_1,\calO_Y),
\calO_Y).$$
This is a reflexive sheaf of rank 1 on $Y$, hence locally free on the smooth locus of $Y$. Let $B$ be the positive Weil
divisor corresponding to $\calL^{-1} = \calL^*$ so that
$$\calL = \calO_Y(-B), \quad \calL^{-1} = \calO_Y(B).$$
We shall call it the {\it branch divisor} of the cover $\pi:X\to Y$.

\med
\proclaim Corollary 1.8. Ouside a finite set of points on $Y$, we have
$$\calH^2(G,\calO_X) = \calO_Y/\Ima(T)\cong \calO_{(p-1)B}.$$

{\sl Proof.} It follows from the definition that $\calL_{p-1} = \Ima T$. Take the open
subset 
$Y'$ where $\calL_1 = \calL$,
then $\calL_{p-1}|Y' = \calO_{Y'}(-(p-1)B') = \Ima(T)|Y'$, where $B' = B|Y'$. This
proves the  assertion.

\med
\proclaim Corollary 1.9. We have
the following formula for the canonical sheaf $\omega_X$ of $X$:
$$\omega_X = \pi^*(\omega_Y\otimes\calL^{-p+1}).$$

{\sl Proof.} Since $\omega_Y$ and $\calL$ are both reflexive sheaves corresponding to some Weil divisors, it is enough
to verify this formula over the complement of a finite set of points in $Y$. By Corollary 1.7, $X$ is a positive divisor in 
$Y\times \bbA^1$ with associated invertible sheaf isomorphic to $p_1^*(\calL^{1-p})$, where $p_1:Y\times \bbA^1\to Y$
is the first projection. Since $\omega_{Y\times\bbA^1} = p_1^*(\omega_Y)$ and $\pi = p_1|X$ the formula 
is just the
adjunction formula. 

\med
{\bf Remark 1.10.} A simple construction of separable $p$-covers $f:X\to Y$
which can be locally given as in Corollary 1.7 over the whole $Y$
(Artin-Schreier covers of simple type in terminology of [Ta]) is given as
follows. Let
$\calL$ be an invertible sheaf on
$Y$ with a non-zero section
$a$ of
$\calL^{1-p}$ and a section $b$ of $\calL^{-p}$. We have the following exact
sequence of sheave in flat topology of $Y$:
$$0\to \alpha_a\to \calL^{-1} \to \calL^{-p}\to 0,$$
where the map $\calL^{-1} \to \calL^{-p}$ is given locally by $x\to
x^p-ax$ and $\alpha_a$ is its kernel. Taking flat cohomology we get an
exact sequence
$$H^0(X,\calL^{-p})\to H_{fl}^1(X,\alpha_a)\to H^1(X,\calL^{-1}).$$
The image of the section $b$ defines an $\alpha_a$-torsor which is
separable $p$-cover ramified outside the divisor of zeroes of $a$.
If $H^1(X,\calL^{-1}) = 0$ every Artin-Schreier cover of simple type can be
obtained in this way. The sheaf $\calL$ is the sheaf $\calL_1$ from
Corollary 1.7.

\bigskip
\centerline {{\bf 2. The case of isolated fixed points.}} 
\med 

In this section only we allow $X$ to be an ''orbitfold K3-surface''.  Recall that
this means that
$\omega_X \cong \calO_X, H^1(X,\calO_X) = 0$  and $\sigma^*(\omega_X) = \omega_{\tilde X}$,
where $\sigma:\tilde X \to X$ is the minimal resolution.
 Let $G = (g)$ be a group of order
$p = \char k$ acting non-trivially on $X$. As before $Y$ denotes the quotient $X/G$
and $\tilde Y$ its minimal resolution.

\med
\proclaim Lemma 2.1. Assume that $G$ acts with finitely many 
isolated fixed
points. Then 
$$\omega_Y \cong \calO_Y.$$

{\sl Proof.} Since we are in characteristic $p >0$ the isomorphism
$\omega_X\cong \calO_X$ is in fact an isomorphism of $G$-sheaves. By
Corollary 1.4, we have 
$$\omega_Y\cong (\pi_*\calO_X)^G \cong \calO_Y.$$ 
\med

Applying the duality theorem on $Y$ we get

\proclaim Corollary 2.2. 
$$\dim_kH^2(Y,\calO_Y) = 1.$$

\med 
\proclaim Lemma 2.3.  The target cohomology
$H^1$ and $H^2$ in the spectral sequences (1.1) and (1.2) are computed as follows:
$$H^1\cong k,$$
$$H^2\cong k^s, s \in \{1,2\}. $$

{\it Proof.}  Recall that every spectral sequence gives the standard five term sequence
$$0 \to E_2^{1,0} \to H^1 \to E_2^{0,1} \to E_2^{2,0} \to H^2.\eqno (2.1)$$
Let us apply it to our second spectral sequence (1.2). Since 
$E_2^{0,1} = H^0(G,H^1(X,\calO_X)) =0$, we get
$$H^1\cong H^{1}(G,H^{0}(X,{\cal O}_X)) = H^1(G,k) = \Hom(G,k) =k.$$
Also, in our situation
$$E_2^{p,q} = H^p(G,H^q(X,\calO_X)) = 0, \quad q\ne 0,2.$$
By [CE], Theorem XV.5.11, we have an exact sequence
$$\ldots E_2^{n,0}\to H^n\to E_2^{n-2,2}\to E_2^{n+1,0}\to H^{n+1}\to
E_2^{n+1-2,2}\to\ldots.\eqno (2.2)$$
Thus the exact sequence (2.1) gives us the exact sequence
$$0\to E_2^{2,0}\to H^2\to E_2^{0,2}\to E_2^{3,0}= k.$$
Since  $G$ acts trivially on $H^{0}(X,{\cal O}_X) = H^{2}(X,{\cal O}_X) = k$, and $T(k)
= 0$, we get 
$$E_2^{2,0} = H^{2}(G,H^{0}(X,{\cal O}_X)) = k/T (k) =k, \quad  
E_2^{0,2}=  H^{0}(G,H^{2}(X,{\cal O}_X))=k.$$ This proves the assertion for $H^2$.

\med

Now we are ready to prove Theorem 1. In fact we prove a more general result since we
are not assuming that $X$ is nonsingular.

\med
\proclaim Theorem 2.4.  Suppose that the fixed locus $X^g$ of $g$
is finite. Then $|X^g|\le 2$ and we have the following possible cases
\item{(i)}$X^g =
\emptyset$: Then
$p = 2$ and
$\tilde Y$ is an Enriques surface and $\tilde X$ is its K3-cover.  
\item{(ii)}$|X^g| = 1$:  $Y$  has one Gorenstein elliptic singularity and $Y$ is a
rational surface, or
$Y$ has one double rational point and $\tilde Y$ is either a K3-surface or an
Enriques surface (the latter case does not happen if $X$ is nonsingular).
\item{(iii)}$|X^g| = 2$: Then $Y$ has two rational double points and
$\tilde Y$ is a K3-surface. 
The cases when $\tilde Y$ is a
K3-surface (resp. Enriques surface) can occur only if
$p
\le 5$ (resp. $p = 2$).

{\sl Proof}. We apply the first spectral sequence (1.1).
By Proposition 1.1 for $i>0$ the sheaf $\calH^i(G,\calO_X)$ is concentrated at a finite set 
of points, so that $E_2^{p,q} = H^p(Y,\calH^q(G,\calO_X)) = 0$ when $q > 0, p > 0$ or 
$p > 2, q = 0$. By [CE],
Proposition XV.5.9, we have the following exact sequence
$$0\to H^{1}(Y,{\cal O}_Y) \to H^1 \to H^{0}(Y,\calH^1(G,\calO_X))\to
H^{2}(Y,{\cal O}_Y)\to $$
$$\to H^2 \to H^{0}(Y,\calH^2(G,\calO_X))\to H^2(Y,\calH^1(G,\calO_X)) = 0
.\eqno (2.3)$$ 
Note first that for $i=1,2$ 
$$H^{0}(Y,\calH^i(G,\calO_X))\cong \oplus_{x \in X^g}H^{i}(G,{\cal O}_{X,x}).$$
As we explained in the proof of Proposition 1.5, for any isolated fixed point $x\in X$
the cohomology
$H^1(G,\calO_{X,x})$ is a non-trivial finite-dimensional vector space over $k$. We also know
that $\dim H^2(Y,\calO_Y) = 1$ and it follows from (2.3) and Lemma 2.3 that 
$\dim H^1(Y,\calO_Y)\le 1.$ 
Let us consider different cases corresponding to all possible values of
$\dim_k H^1(Y,\calO_Y)$.

\medskip\noindent
Case 1: $\dim_k H^1(Y,\calO_Y) = 1$.

Using (2.3), we get $H^0(Y,\calH^1(G,\calO_X)) = 0$ or $k$ so that $G$ has either acts 
freely or has one fixed point.

\medskip\noindent
Case 2: $\dim_k H^1(Y,\calO_Y) = 0$.

We get $H^0(Y,\calH^1(G,\calO_X)) = k$ or $k^2$. In this case, $X^g$ consists of one
or  two points.

Now let us investigate each of the two cases in more detail.

\medskip
{\bf Case 1}: If $H^0(Y,\calH^1(G,\calO_X)) = 0$, the map $\pi:X\to Y$ is \'etale, and
hence, the lift of the action to $\tilde X$ is free and $\tilde Y =\tilde X/G$
 is an Enriques surface. Since $K_{\tilde Y} = 0$, this can happen only when $p = 2$.
(see [CD], Theorem 1.1.3). Since the K3-cover is separable $\tilde Y$ is a
non-classical Enriques surface of 
$\mu_2$-type and
$p = 2$ (see loc. cit., p.77).

If $H^0(Y,\calH^1(G,\calO_X)) = k$, then the image $y =\pi(x)$ of the unique fixed point 
$x\in X^g$  is a rational singularity ([Pe1], Theorem 6).
 Let $V = Y\setminus \{y\}$ and $U = \pi^{-1}(V)$. Since $\pi^*(\omega_V) = \omega_U =
{\cal O}_U$ and the homomorphism $\Pic(V)\to \Pic(U)$ is injective (its kernel is
isomorphic to $H^1(G,k^*) = 0$), we see that $\omega_V = {\cal O}_V$. Since $Y$ is
normal, this implies that $\omega_Y = \calO_Y$. So, $y$ is a  rational double point.
This implies that the canonical class of 
$\tilde Y$ is trivial and $H^1(\tilde Y,\calO_{\tilde Y}) = 0$. By the
classification of algebraic surfaces, $\tilde Y$ must be a non-classical Enriques
surface. This could happen only if $p = 2$. Let us show that this is impossible if $X$
is nonsingular.

 Let $\tilde X\to \tilde Y$ be the K3-cover of $\tilde Y$. It is a principal
bundle with respect to a group scheme ${\cal A}$ of order 2 isomorphic to either
$\bbZ/2$ or
$\alpha_2$ (see [CD], Chapter 1, \S3).  Restricting the cover
over the complement of the exceptional locus $E\subset \tilde Y$, we get a principal
${\cal A}$-cover of $Y\setminus\{y\}$. Its pre-image under the map $p:X\setminus
\{x\}\to Y\setminus \{y\}$ is a principal ${\cal A}$-cover of $X\setminus \{x\}$. Since
$x$ is a nonsingular point, we can apply the purity theorem ([CD], 0.1.10) to
extend this cover to a principal ${\cal A}$-cover of $X$. Since $H^1(X,{\cal O}_X)
= 0$ this cover must be trivial (use loc. cit, Proposition 0.1.7 and 0.1.9).
This implies that $X\setminus\{x\} \cong \tilde X\setminus E'$, where $E'$ is the
pre-image of $E$ in $\tilde X$. Now, if $\calA = \bbZ/2$, $\tilde X$ is a nonsingular
K3-surface and
$E'$ is the disjoint union of two curves. Resolving the points of indeterminacy
of the rational map
$\tilde X\to X$ by blowing up some points on
$E'$, we get a birational morphism $\tilde X'\to  X$ which has disconnected
pre-image over $x$. This contradicts  Zariski's Connectedness Theorem. Next we assume
that ${\cal A} =\alpha_2$. In this case $X\setminus E'$ is an inseparable cover of
$\tilde Y\setminus E$. Hence the $l$-adic Euler characteristic with compact
support $e_c^{l}$ of
$\tilde X\setminus E'$ is equal to 
$$e_c^{l}(\tilde Y\setminus E) = e_c^{l}(\tilde Y)-e_c^{l}(E) = 12-e_c^{l}(E) \le
10.$$
On the other hand, $e_c^{l}(X\setminus \{x\}) = e_c^{l}(X)-1 = 23.$ This contradiction
proves the claim.

\medskip
{\bf Case 2}: Arguing as in the previous case we have
the following possible cases:
\item {(a)} $\tilde Y$ is a rational surface and $Y$
has one non-rational Gorenstein singularity. 
An easy argument using the Leray spectral sequence for the minimal resolution $\tilde
Y\to Y$ shows that this is an elliptic singularity.
\item{(b)} $\tilde Y$ is a K3-surface and $Y$ has one rational double point
with
$H^0(Y,\calH^1(G,\calO_X))= k^2$
\item{(c)} $\tilde Y$ is a K3-surface and $Y$ has two rational double points with
$H^0(Y,\calH^1(G,\calO_X)) = k^2$.
\item{(d)} $\tilde Y$ is a K3-surface and $Y$ has one rational double point with
$H^0(Y,\calH^1(G,\calO_X)) = k$.

 Now, each of the rational double points on
$Y$ has a smooth covering of order $p$ \'etale above a punctured neighborhood of the point. 
By Artin ([Ar2], Corollary 2.7), this is possible only if $p\le 5$.
This finishes the proof.

\med
{\bf Remark 2.5.} If $p\ge 3$, then the nonemptyness of $X^g$ can be seen
more easily;  the algebraic Euler characteristic $\chi(X,\calO_X)$ is equal to
$2$, which is not divisible by
$p$.

\med
{\bf Remark 2.6.} In the case $p = 2$ Artin shows in
[Ar1] that the completion of the local ring of the image $y\in Y$  of an
isolated fixed point $x\in X$ is isomorphic to 
the ring $k[[x,y,z]]/(z^2+abz+a^2y+b^2x)$, where $a,b\in k[[x,y]]$ are
relatively prime nonunits. Also it follows that  the image of the trace
map
$\Ima(T)$ equals the ideal $(a,b,z)$. Thus, by Proposition 1.5,
$$\dim H^1(G,\calO_{X,x}) = \dim H^2(G,\calO_{X,x}) = \dim
\calO_{Y,y}/(a,b,z).$$
This shows that in cases 2(c),(d) the completion of the local ring
$\calO_{Y,y}$ is isomorphic to the ring
$k[[x,y,z]]/(z^2+xyz+x^2y+y^2x)$. This is a rational double point of
type $D_{4}^{(1)}$ from Artin's list in [Ar2]. In case 2(b), we get
the ring $k[[x,y,z]]/(z^2+x^2yz+x^4y+y^2x),$ or
$k[[x,y,z]]/(z^2+x^2yz+x^5+y^3)$. This is a  double rational point of
type $D_{8}^{(2)}$ or $E_8^{(4)}$, respectively. 

The case 2(a) does not occur in characteristic 2. However it may occur in
characteristic 3 as the following example from [Pe1] shows. The group
$G =\bbZ/3\bbZ$ acts on $k[[u,v]]$ by $\sigma(u) = u+y^i, \sigma(v) = v+u$.
The ring of invariants is isomorphic to the ring
$k[[x,y,z]]/(z^3+y^{2i}z^2-y^{3i+1}-x^2).$ The ideal $\Ima(T)$ equals 
$(x,y^i,z)$. When $i = 1$ we get a rational double point of type 
$E_6^{(1)}$. When $i = 2$ we get an elliptic singularity with $\dim
H^2(G,\calO_{X,x}) = 2$. 

Finally note that a rational singularity of $Y$ need not be a double point.
This may occur already when $p = 3$ (see Example 13 from [Pe1]).

\med
{\bf Remark 2.7.} We remark that even if $|X^g|=1$, or $2$, the {\it fixed point
scheme} $X^g$  can not be isomorphic to $\Spec(k)$, locally. More generally, 
if a finite abelian $p$-group acts on a nonsingular
variety of positive dimension defined over $k$ of characteristic $p$ by 
$k$-automorphisms, then its fixed point scheme can not
contain any isolated closed point (see  [ABK],
Theorem 3.1). This can
 be seen by noting that at an isolated fixed point the linear 
action of the group on the 
tangent space of the fixed point is unipotent, and hence fixes a nonzero tangent vector as well.

\med
{\bf Examples 2.8.} 
1. Let $X \subset (\bbP^1)^4$ be a complete intersection of two hypersurfaces 
$(F_1 =0)$ and $(F_2 =0)$ of multidegree $(1,1,1,1)$. Let $g$ be the involution
on $(\bbP^1)^4$,
$g(x,y,z,w)=(y,x,w,z)$. Suppose that $g^{*}F_i = F_i$, $i=1,2$. Then 
$|X^g| =8$ 
if $\char k=0$ and  $|X^g| =2$ if $\char k=2$.
\med
2. Let $X \subset \bbP^4$ be a complete intersection of a quadric and a cubic, both
invariant under the action 
$g(x_0,x_1,x_2,x_3,x_4)=(x_0,x_2,x_1,x_4,x_3)$ of order 2. Suppose that the cubic
has the term $x_0^3$. Then $|X^g| =1$ if $\char k=2$.
\med
3. Let $X$ be the Fermat quartic $x_0^4+x_1^4+x_2^4+x_3^4=0$ in $\bbP^3$. Let $g$ be
the automorphism of order $3$, $g(x_0,x_1,x_2,x_3)=(x_0,x_2,x_3,x_1)$. Then  $|X^g| =6$
if $\char k=0$ and  $|X^g| =1$ if $\char k=3$.
\med
4. Let $X$ be the hypersurface in $(\bbP^1)^3$ of multidegree $(2,2,2)$, invariant 
under the action of $g:(x,y,z)\to(y,z,x)$. Then $|X^g| =6$
if $\char k=0$ and  $|X^g| =2$ if $\char k=3$.
\med
5. Let $X \subset \bbP^5$ be a complete intersection of three quadrics,
invariant under the action 
$$g(x_0,x_1,x_2,x_3,x_4,x_5)
    = (x_0,x_2,x_3,x_4,x_5,x_1)$$ 
of order 5. Suppose that one of the three quadrics has the term  $x_0^2$.
Then $X^g =\{(0,1,1,1,1,1)\}$ if $\char k=5$.
\med
6. Let $X \subset \bbP^6$ be a complete intersection of three quadrics and 
a hyperplane, invariant under the action 
of order 7, $g(x_0,x_1,x_2,x_3,x_4,x_5,x_6)=(x_1,x_2,x_3,x_4,x_5,x_6,x_0).$ 
Then $X^g =\{(1,1,1,1,1,1,1)\}$ if $\char k=7$.
\bigskip

\centerline{{\bf 3. The case of non-isolated fixed points.}} 
\med
From now on we assume that the set $X^g$ of fixed points of the
involution 
$g$ contains a one-dimensional part. Denote it by $F$. Recall from Section $1$
that we have defined the branch divisor $B$. Its support
$B_{red}$ is equal, outside a finite set of points, to the support of the sheaf
$\calH^2(G,\calO_X)$.

\med
\proclaim Lemma 3.1.
$$B_{red} = \pi(F).$$

{\sl Proof.} Let $U$ be the open subset where $B_{red}$ is equal to the support of 
$\calH^2(G,\calO_X)$. 
Over $U$, we have by  Corollary 1.2, $B_{red}\subset \pi(X^g)$. Let $x\in X^g$, then 
$g$ acts on the local ring $\calO_{X,x}$ sending its maximal ideal ${\frak
m}$ to itself. Since $T$ acts trivially on constants, the image of $T$ on
$\calO_{X,x}$ is contained in ${\frak m}\cap
\calO_{Y,\pi(x)}$. Thus
$\pi(x)$ belongs to $B_{red}$. Now two one-dimensional closed subsets $B_{red}$ and $\pi(F)$ are 
equal on $U$, hence are equal over
the whole $Y$ since $Y\setminus U$ is of codimension $\ge 2$. 

\med
\proclaim Proposition 3.2. Let $\omega_Y$ be the dualizing sheaf
 of
$Y$. Then
$$\omega_Y = \calO_Y(-(p-1)B).$$
Moreover
$$\calO_{(p-1)B}\cong \calO_Z,$$
where $Z$ is the closed subscheme of $Y$ defined by the ideal $\Ima(T)$.

 {\sl Proof.} As we saw in the proof of
Lemma 2.1,
$(\pi_*\omega_X)^G=\calO_Y$ as $G-\calO_Y$-modules.  Applying Proposition 1.3,
we obtain the exact sequence
$$0\to \omega_Y\to \calO_Y\to {\cal E}xt_{\calO_Y}^1(\calH^2(G,\calO_X),\omega_Y)\to
0.$$
Let $Z$ be the closed subscheme of $Y$ defined by the ideal sheaf
$\Ima(T)$. Then $\calH^2(G,\calO_X) = \calO_Z$ and, by the duality
${\cal E}xt_{\calO_Y}^1(\calH^2(G,\calO_X),\omega_Y) \cong \omega_Z$.
Let $W$ be the closed subscheme of $Y$ defined by the ideal
sheaf $\omega_Y\subset \calO_Y$. Restricting to the smooth locus
$Y'$ of
$Y$ we find that
$Z' = Z|Y'$ and $W'= W|Y'$ are Cartier divisors and hence 
$$\calO_{W'} \cong \omega_{Z'} = \omega_{Y'}(Z')/\omega_{Y'} \cong
\omega_{Y'}(Z')\otimes \calO_{Z'} .$$
Since $\omega_{Y'}(Z')$ is an invertible sheaf on $Y'$, we see that the
subschemes $Z'$ and $W'$ are identical. The corresponding Cartier
divisor extends uniquely to a Weil divisor $Z= W$ on $Y$. Since $(p-1)B = Z$ on a set 
with complement of codimension $\ge 2$ we see that $Z = (p-1)B$.

\proclaim Corollary 3.3. Let $K_Y$ be the canonical divisor of $Y$. Then
$$K_Y\sim -(p-1)B.$$

{\sl Proof.} By Corollary 1.9, we have 
$$\pi^*(\omega_Y\otimes\calO_Y((p-1)B)) \cong \calO_X.$$
Since $\pi^*:\Pic(Y)\to\Pic(X)$ is injective, we get the desired formula.    
\med

\proclaim Proposition 3.4. We have
$$H^1(Y,\calO_Y) = 0.$$

{\sl Proof.}  Consider the exact sequence
$$ 0\to \calO_Y\to \pi_*(\calO_X)\to \calE\to 0.\eqno (3.1)$$
Since
$\pi:X\to Y$ is a finite morphism, we have
$H^1(X,\calO_X) = H^1(Y,\pi_*(\calO_X)) = 0,$
and hence we infer from (3.1) that $H^1(Y,\calO_Y) = 0$ if $H^0(Y,\calE)
= 0$. Aplying Proposition 1.6 we find that $\calE$ admits a filtration with quotients 
isomorphic to ideal sheaves $\calL_i$
and $\calL_i\subset \calL_1$. So, it is enough to show that $H^0(Y,\calL_1)
= 0$. But this is easy. Over a complement to a finite set of points,
$\calL_1 = \calO_Y(-B)$ for some positive divisor $B$. Obviously such a
sheaf has no nonzero sections.

\proclaim Proposition 3.5. For any $i > 0$ we have 
$$\dim H^0(Y,\calH^i(G,\calO_X)) = \dim H^1(Y,\calH^i(G,\calO_X)) = 1.$$

{\sl Proof.} Applying the exact sequence (2.1) and Lemma 2.3 (both do not
use the assumption that $G$ acts with finitely many fixed points) we find
that $ \dim H^0(Y,\calH^1(G,\calO_X)) = 1$. Since, for any $q\ge
0, p > 1$,
$'E_2^{p,q} = H^p(Y,\calH^q(G,\calO_X)) = 0$, we can apply Proposition XV.5.5
from [CE] to find an exact sequence
$$0\to E_\infty^{1,1}\to H^2\to E_\infty^{0,2}\to 0.$$
It is easy to see that  $E_\infty^{p,q} = E_2^{p,q}$ for $(p,q) =
(1,1), (0,2)$. This gives us the exact sequence
$$0\to H^1(Y,\calH^1(G,\calO_X))\to H^2\to H^0(Y,\calH^2(G,\calO_X))\to
0.\eqno (3.2)$$ 
Applying Proposition 3.2 and duality, we get
$$\chi(Y,\calH^2(G,\calO_X)) = \chi(\calO_Z) = -\chi(\omega_Z) = \chi(\calO_{(p-1)B}) =$$
$$\chi(Y,\calO_Y)-\chi(Y,\calO_{Y}(-(p-1)B)) = \chi(Y,\calO_Y)-\chi(Y,\omega_Y) = 0.$$
Thus, applying Proposition 1.5, we get
$\dim H^1(Y,\calH^1(G,\calO_X)) = 1$, and the exact sequence (3.2) together
with Lemma 2.3 proves the assertion.
\med

\proclaim Corollary 3.6. The fixed locus $X^g$ is connected  
and its image in $Y$ is equal to the support of the branch divisor $B\in |-{1\over p-1}K_Y|$ 
with $$\dim H^i((p-1)B,\calO_{(p-1)B}) = 1, i  = 0,1.$$

\med
Let us summarize what we know about the quotient surface $Y$. 

\proclaim Theorem 3.7. Assume that $G$ acts on $X$ with non-isolated fixed points. Then 
the quotient surface $Y = X/G$ is a normal  surface with effective  anticanonical divisor. 
It is rational and has at most rational singularities.

{\sl Proof.} Let $\sigma:\tilde Y \to Y$ be a resolution of singularities and
$\tilde X$ be the normalization of $\tilde Y$ in the field of rational
functions of $X$. Let $Z\to \tilde X$ be its resolution of singularities.
Then the composition
$f:Z\to \tilde X\to X$ is a resolution of singularities of a nonsingular
surface. This easily implies that $R^1f_*\calO_Z = 0$ and this gives
immediately that $H^1(Z,\calO_Z)\cong H^1(X,\calO_X) = H^1(\tilde
X,\calO_{\tilde X}) = 0.$ Now the group
$G$ acts on $\tilde X$ with the quotient equal to $\tilde Y$. One can
extend the proof of Proposition 3.4 to this situation to prove that 
$H^1(\tilde Y,\calO_{\tilde Y}) = 0.$ This means that singularities of $Y$
are rational. Also it is known that 
$$K_{\tilde Y} = \sigma^*(K_Y)-\Delta\eqno (3.3)$$
for some positive divisor $\Delta$ supported on the exceptional locus of
$\sigma$. Thus $-K_{\tilde Y}$ is effective, and $\tilde Y$ is of Kodaira
dimension equal to $-\infty$. Since $H^1(\tilde Y,\calO_{\tilde Y}) = 0$,
$\tilde Y$ is a nonsingular rational surface.

\med
Recall that a normal surface $S$ is called {\it anticanonical} if $-K_S$ is
effective. Thus our quotient surface $Y$ is a normal rational anticanonical
surface. It follows from the formula (3.3) that its minimal resolution is a
nonsingular rational anticanonical surface. Now we can invoke the
classification of such surfaces (see [Sa1, Sa2]). First, they are divided
into the classes corresponding to different anti-Kodaira dimension
$\kappa^{-1}(S)$, i.e. the
$D$-dimension of $Y$, where $D = -K_S$. For any anticanonical surface,
$\kappa^{-1}(S)\ge 0$, so we have three possible cases for our $Y$:
$\kappa^{-1}(Y) = 0,1,2$.

Next, we use the Zariski decomposition 
$$-K_S = P+N,\eqno (3.4)$$
where $P$ is a nef $\bbQ$-divisor, and $N$ is an effective $\bbQ$-divisor such
that either $N = 0$ or its irreducible components $N_i$ define a negative
definite intersection matrix and $P\cdot N_i = 0$ for each $i$. Here we use the
standard intersection theory on a normal surface due to Mumford. We set
$$\nu(S,-K_S) = \cases{0&if $P\equiv 0$\cr
1&if $P^2 = 0, P\not\equiv 0$\cr
2&if $P^2 > 0$\cr},$$
$$d(S) = P^2.$$
We have
$$d(S) =\cases{0 &if $\kappa^{-1}(S) = 0,1$\cr
> 0&if $\kappa^{-1}(S) = 2$\cr}.$$
By contracting irreducible curves $C$
with
$C^2 < 0, C\cdot K_S < 0$ we arrive at a minimal normal anticanonical
surface
$S_0$. Also by resolving rational double points (which do not affect the
canonical class) we may assume that either $S_0$ is nonsingular or contains at
least one singularity which is not a rational double point.
\med
{\bf Remark 3.8.} The occurence of singularities of the quotient surface
depends on the linear part of the action of $g$. For example, if $p = 3$ and
the linear part of $g$ in a neighborhhod of a point $x\in F$ has one Jordan 
block then the image of $x$ is a nonsingular point (see Corollary 5.15 in
[Pe2]).

\med
\proclaim Theorem 3.9 ([Sa1]). Let $S$ be a nonsingular anticanonical
rational surface. There exists a birational morphism
$\phi:S\to S_0$ where $S_0$ is a nonsingular anticanonical surface with $\nu(S)
= \nu(S_0,-K_{S_0})$ and Zariski decomposition $-K_{S_0} = P_0+N_0$ where
$n_0\phi^*(P_0) = nP$ for some positive integers $n$ and $n_0$. We have four
possibilities:
\item{(i)} $\nu(S,-K_S) = \kappa^{-1}(S) = 0$:$P_0 = 0$;
\item{(ii)}$\nu(S,-K_S) = 1,\kappa^{-1}(S) = 0$: $P_0 = -C_0$ where $C_0$ is an indecomposable
curve of canonical type (i.e. each its irreducible component $R$ satisfies
$R\cdot K_{S_0} = C_0\cdot R = 0$ and $C_0$ is not a sum of two curves
satisfying such property) with normal sheaf
$\calN_{C_0}$ of infinite order in
$\Pic(C_0)$;
\item{(iii)}$\nu(S,-K_S ) = \kappa^{-1}(S) = 1$: there exists a minimal elliptic (or
quasi-elliptic) fibration on
$S_0$ for which
$-mK_{S_0}$ is linearly equivalent to a fibre for some $m > 0$.
\item{(iv)} $\nu(S,-K_S) = \kappa^{-1}(S) = 2$.

\med
Note that in our situation the surface $Y$ is often singular so there is
no hope for a complete classification and even when singularities are
quotient singularities this is sill hopeless (see [Ni3]). In the case when
$Y$ is Gorenstein (i.e. $\omega_Y$ is locally free), and we are in case
(iv) , it is obtained from a Del Pezzo surface $V$ by blowing down
$(-2)$-curves. Notice
$-K_V$ is divisible by $p-1 > 1$ in $\Pic(V)$ only when $p = 3$ and $V = Y$
is a quadric. 
\med
  Let us return to our situation. We have

\med
\proclaim Lemma 3.10. Let $D$ be a nef divisor on a K3-surface. If $|D|$ has
a fixed part, then there exists an irreducible curve $E$ of arithmetic
genus 1 such that $D\cdot E = 1$ and
$D = aE+\Gamma$, where $a \ge 2$ and $|E|$ is a free pencil and $\Gamma$ is a smooth 
rational curve with $E\cdot
\Gamma = 1$. If $|D|$ has no fixed part and $D^2 = 0$, then $D = aE$ for some free pencil 
$|E|$. Moreover, if $D^2 > 0$ and $|D|$ has no fixed part, then $|D|$ has no base-points 
and the corresponding map 
$\phi_{|D|}$ is either of degree 2 or birational onto a surface with at most
double rational points as singularities.

{\sl Proof.} This is a well-known result of B. Saint-Donat 
(see, for example, [Re], 3.8, 3.15 ).
\med

\proclaim Theorem 3.11. 
$$\nu(Y,-K_Y) = \kappa^{-1}(Y) = \kappa(X,X^g).$$ 

{\sl Proof.} We  put $F = X^g$. Let $B$ be the branch divisor of $\pi:X \to Y$
so that $aF =\pi^{*}(B)$ for some $a > 0$. From the projection formula we find
$$\calO_Y(nB) \subset \pi_* \pi^*(\calO_Y(nB)) = \calO_Y(nB)\otimes \pi_*(\calO_X).$$
Thus 
$$\dim H^0(X,\calO_X(naF)) = \dim H^0(Y,\pi_* \pi^*(\calO_Y(nB))) = 
\dim H^0(Y,\calO_Y(nB)
\otimes \pi_*(\calO_X)).$$
By Corollary 1.7, outside a finite set of points, $\pi_*(\calO_X)$ admits a
filtration with successive quotients equal to
$\calL^{i}, i = 0,\ldots,p-1,$ for some invertible sheaf $\calL$. Also we know that 
$\calL \cong \calO_Y(-B).$ This immediately implies that
$$\dim H^0(Y,\calO_Y(nB)) \le \dim H^0(X,\calO_X(naF)) \le p\dim H^0(Y,\calO_Y(nB)).$$   
In particular, 
$$\kappa^{-1}(Y) = \kappa(Y,B) = \kappa(X,F).$$
Let $F = P+N$ be the Zariski decomposition of $F$. Let $m$ be an integer such that 
$mP\in \Pic(X)$. Then
$mF = mP+mN$ where $mP$ is a nef divisor. Obviously, $mN$ is the negative part of the Zariski
decomposition of $mF$. Thus, $mN$ is contained in the fixed part of $|mF|$. 
Assume $P^2 = 0$, i.e. $\nu(X,F) = 1$. By Lemma 3.10,
$mP$ has no fixed part and is composed of a pencil. Thus $mP$ is the moving part of $mF$ 
and $\kappa(X,F)
= \kappa(X,mP) = 1$. Assume $P^2 > 0$. It follows from Lemma 3.10 that
$|2mP|$ has no fixed part. So,
$\nu(X,F) = \kappa(X,F) =\kappa(X,mP) =2$.
 On the other
hand, by Lemma 2.5 of [Sa1] $P = \pi^*(P'), N = \pi^*(N')$ where $B = P'+N'$ is the Zariski
decomposition of $B$. This implies that $\nu(Y,B) = \nu(X,F)$ and finishes the proof.

\med
So, we see that there are three different cases to consider according to three possible 
values of the Kodaira dimension of the divisor $F$ of fixed points of $G$. 
This will be dealt with in the next sections.

\bigskip
\centerline{{\bf 4. $\kappa(X,F) = 0$}}
\med
In this case the Zariski decomposition of $F$ consists only of the negative part $N$. 
Each component of
$F$ must be a smooth rational curve, and the sublattice of $\Pic(X)$ spanned by the
classes  of the components is a lattice of type A-D-E. It is known that we can blow
down $F$ to  a double rational singular point $x$ of an ``orbitfold K3-surface '' $X'$.
The involution $g$ acts on $X'$ with one isolated
fixed point $x$. Then we apply the results of section 2.  By Theorem 2.4, we see that
the  quotient $Y' =X'/G$ has either one rational singular point or one elliptic
singular point. 
  A rational double point has a rational double
point as  a covering of order $p$ which is \'etale above a punctured neighborhood of 
the point, 
only if $p\le 3$ ( see [Ar2]). Actually, Artin shows that $x$  must be of
type $A_{8r-4n-1}$,
$A_{8r-4n+1}$,
$A_{2}$, $D_{4}^{(1)}$ if $p=2$, $A_{1}$, $D_{4}$ if $p=3$, and its
corresponding  rational double point in $Y'$ is of type $D_{2n}^{(r)}(2r\ge
n)$, $D_{2n+1}^{(r)}(2r\ge n)$,
$E_6^{(1)}$, $E_7^{(3)}$ if $p=2$,  $E_7^{(1)}$, 
$E_8^{(2)}$ if $p=3$, respectively. 

This proves Theorem 2 from the introduction.

\med
{\bf Example 4.1} Consider the K3-surface $X$ over a field of characteristic
2 given by the equations $$x_0^3+x_1^3+x_2^3+x_3^3+x_4^3 = 0,\quad x_0^2+x_1x_2+x_3x_4 = 0.$$
It has the involution defined by $(x_0,x_1,x_2,x_3,x_4)\to
(x_0,x_2,x_1,x_4,x_3)$ and has the unique fixed point
$(0,1,1,1,1)$ which happens to be an ordinary double point of $X$. The
quotient surface $Y$ is isomorphic to a surface of degree 3 in the weighted
projective space $\bbP(1,1,1,2)$ given by the equation
$$x_0^3+x_1^3+x_2^3+x_1x_3+x_2(x_0^2+x_3)= 0.$$
Recall that $\bbP(1,1,1,2)$ is isomorphic to the projective cone over the
Veronese surface in $\bbP^5$. Our surface $Y$ passes through the vertex of
this cone. It is a rational surface with an elliptic singularity. Note that
by adjunction formula $\omega_Y = \calO_{\bbP}(-2)\otimes\calO_Y$. There is
no contradiction here with Lemma 2.1 since $\calO_{\bbP}(-2)$ is not an
invertible sheaf and $\calO_{\bbP}(-2)\otimes\calO_Y\ne \calO_Y(-2)$.

\bigskip
\centerline{{\bf 5. $\kappa(X,F) = 1$}}

\med

\proclaim Lemma 5.1. There exists a divisor $F'$ with the same support as $F$ such that
the linear system
$|F'|$ defines an elliptic or quasi-elliptic fibration $\phi:X \to \bbP^1$.

{\sl Proof.} Let $F= P+N$ be the Zariski decomposition of $F$. As in the proof of
Theorem 3.11 we write the Zariski decomposition of integral divisors 
$mF = P'+N'$, where
$P'$ is a nef divisor and
$N'$ is a sum of smooth rational curves whose intersection matrix is negative
definite, and 
hence of type A-D-E. Since $F$ is connected, all components of $N'$ are contained in
$P'$. 
In particular the support
of $P'$ is equal to $F$. It follows from Lemma 3.10 that 
$P'^2 = 0$ and  $|P'| = |aE|$ for some pencil $|E|$ of curves of arithmetic
genus 1. Since $P'$ and $F$ have the same support we are done.
\med

\proclaim Lemma 5.2. The induced map $g^*$ on the base curve $\bbP^1$ of
the fibration $\phi:X \to \bbP^1$ is of order $p$ unless $p = 2$ and the fibration is
elliptic.

{\sl Proof.} Assume first that $\phi:X \to \bbP^1$ is a 
quasi-elliptic fibration. 
 In this case a general fibre is an irreducible curve of arithmetic genus $1$
with a cusp. The closure of the cusps of irreducible fibres is a curve $C$, called the
cuspidal curve. It is known that the
restriction of $\phi$ to
$C$ is a purely inseparable cover of degree $p$. Thus $C\cdot E = p > 1$. 
If fibres were preserved under $g$, the curve $C$ would be contained in $X^g$. However
by the previous lemma this is impossible. 

 Now suppose that $\phi:X \to \bbP^1$ is an
elliptic fibration and that $g^*$ acts as identity on the base curve $\bbP^1$. 
Then $g$ becomes an automorphism of the elliptic  curve
$X/\bbP^1$ over the function field of  $\bbP^1$. 
The automorphism $g$ induces  an automorphism $\bar{g}$ of order $p$ on
the Jacobian of this elliptic  curve. It
fixes a specific fibre pointwisely. Let $j:J\to\bbP^1$ be the Jacobian
surface. Note that $J$ is a K3-surface (cf. [CD], Theorem 5.7.2).  By
Lemma 5.1,
$\bar{g}$ can not fix a section of $j$ pointwisely, and so is a translation by a
$p$-torsion section $P$.
We recall the height paring $<-,->$ defined by Shioda
on the Mordell-Weil group of an elliptic surface. Since $P$ is a torsion section, 
$<P,P>=0$. On the other hand, explicit formula for the height paring
(Shioda [Sh], Theorem 8.6) tells that in our case
$$<P,P> = 4+2P\cdot O - \sum_{v}{\rm contr}_{v}(P),$$
where $O$ is the zero-section of $j:J\to\bbP^1$, $P\cdot O$ the intersection number 
of the divisors $P$ and $O$, the summation over all critical values of $j$, and 
${\rm contr}_{v}(P)$
is a rational number determined by the incidence relation among the divisor $P$ and 
irreducible components of the singular fibre $j^{-1}(v)$, e.g. it is 
equal to $0$ if $P$ meets the identity component of  $j^{-1}(v)$. Assume for
simplicity that
$p\ge5$. Since $P$ is a 
torsion section of order $p\ge5$, the number ${\rm contr}_{v}(P)$ is not zero only if
the singular fibre $j^{-1}(v)$ is of type $I_{ap}$. Let $I_{a_{1}p}$,..., $I_{a_{r}p}$
be
all singular fibres of $j$ of such type. Then by the formula for ${\rm contr}_{v}(P)$
([Sh], p. 229 or  [CZ]), we have
$$ \qquad \sum_{v=1}^{r}(i_{v}(p-i_{v})a_{v}/p) = 4 + 2P\cdot O,\eqno
(5.1)$$
 where we assume that $P$ meets the $i_{v}a_{v}$-th component of
$j^{-1}(v)$ counted  from the
identity component. From the bounds for the second Chern number and Picard number of
the 
surface $J$, we have
$$24\ge \sum_{v=1}^{r}pa_{v} \quad {\rm and} \quad 20\ge
\sum_{v=1}^{r}(pa_{v}-1).\eqno (5.2)$$
Under these constraints, some computation shows that the left hand side of (5.1)
takes its 
maximum $4$. On the other hand, we must have $P\cdot O \ge1$, because the section $P$ 
intersects the zero-section $O$ at the fibre fixed pointwisely by $\bar{g}$. This is a
contradiction to (5.1), and we have proved the assertion for $p\ge5$.
If $p = 3$, the number ${\rm contr}_{v}(P)$ is not zero only if
the singular fibre $j^{-1}(v)$ is of type $I_{ap}$, $IV$ or $IV^*$, and hence (5.1)
and (5.2)
can be written as follows:
$$ \qquad 2s/3 + 4t/3 + \sum_{v=1}^{r}(i_{v}(3-i_{v})a_{v}/3) = 4 + 2P\cdot O,$$
$$24\ge 4s + 8t +3\sum_{v=1}^{r}a_{v} \quad {\rm and} 
\quad 20\ge 2s+6t+\sum_{v=1}^{r}(3a_{v}-1),$$ 
where $s$(resp. $t$) is the number of fibres of $j$ of type $IV$(resp. $IV^*$). 
A similar computation also leads to a contradiction in this case.

\med
Next example shows that the assumption  $p\ne 2$ is essential in the case of elliptic
fibrations.

\med
{\bf Example 5.3.} Here $p = 2$. Consider a rational elliptic surface
$f:V\to \bbP^1$ from the list of extremal rational surface in [La] which is
given by the Weierstrass equation
$y^2+txy+ty=x^3$. It has three degenerate fibres
of type $A_5,A_1$ and
$A_2^*$. Let us apply the construction from Remark 1.10 by taking $\calL =
f^*(\calO_{\bbP^1}(-1))$. Take its section $a$ defined by the fibre $V_0$ of
type
$A_2^*$ and take $b$ to be the sum of two disjoint nonsingular fibres.
Let $X'\to V$ be the corresponding cover. By the formula for the canonical
class (Corollary 1.9) we get that $\omega_{X'} = 0$. It has an elliptic
fibration and unramified outside $V_0$. The preimage of the other singular
fibres are fibres of type $A_5,A_5,A_1,A_1$. They satisfy the condition 
(5.1). One
can show by local computations using equations from Corollary 1.7 that
$X'$  has a  rational double point of type $D_4^{(1)}$ over the singular
point of
$V_0$. After its resolution we get an elliptic K3-surface $\phi:X\to \bbP^1$
with 4 fibres of multiplicative type as above and one fibre of additive type
$E_6$. The Mordell-Weil group of $f$ is isomorphic to $\bbZ/6\bbZ$ (see for example
[CD], Corollary 5.6.7). This shows that the 2-torsion part of the Mordell-Weil group
of
$\phi$ is is not trivial. The corresponding translation automorphism has fixed point
set equal to the fibre of type $E_6$ and preserves the fibres of the elliptic
fibration.

\proclaim Lemma 5.4. Assume $p\ne 2$. The genus 1 fibration $\phi:X \to
\bbP^1$ has at least
$p$ degenerate fibres away from the fibre containing $F$.

{\it Proof.} Since $g$ acts on the base with one fixed point, we obtain that other
singular fibres form the union of orbits of $G$ each consisting of $p$ fibres. 
So it is enough to show that $\phi$ contains more than one degenerate fibre. Assume
first that $\phi$ is quasi-elliptic. Then the formula for the Euler characteristic of
$X$ from [CD], Proposition 5.1.6 gives that the unique degenerate fibre must have 21
irreducible components. The only type of fibre of additive type with so many 
components
is of type $D_{20}$ (Kodaira's $I_{16}^*)$. However such a fibre may occur in a
quasi-elliptic fibration only if $p = 2$(loc. cit., Corollary 5.2.4). Assume now that
$\phi$ is an elliptic fibration with one degenerate fibre. Over the complement to one
point of the base, 
$\phi$ defines an abelian group scheme $\calA'\to \bbA^1$. For any prime
$l\ne p$ its group of $l$-torsion sections defines an unramified extension of 
the affine line $\bbA^1$ of degree $l^2$. 
Since the algebraic fundamental group of the affine line in characteristic $p > 0$ is a
profinite $p$-group, we obtain that all $l$-torsion points are defined over
$k(\bbP^1)$. Thus $\phi$ has sections of any order $l\ne p$ which contradicts the
Mordell-Weil Theorem.

\med

\proclaim Lemma 5.5. Let $\tilde Y$ be a minimal resolution of the quotient
surface $Y = X/<g>$. Then $\tilde Y$ is a rational elliptic surface.

{\sl Proof.} We already know that $Y$ is a rational surface (Theorem 3.7). 
Away from the branch divisor $B$, it has an elliptic fibration induced by 
the elliptic fibration $\phi:X\to \bbP^1$.
Let $\sigma:\tilde Y\to Y$ be a minimal resolution of $Y$. Then it is a
rational elliptic surface which however may not be (relatively) minimal. 
By blowing down the exceptional components of the fibre
$\bar E$ arising from the fixed locus $F$ of $X$, 
we obtain a minimal rational elliptic surface $V$. 
So, we are exactly in the situation described by Theorem 3.9, case (iii).

\med
\proclaim Corollary 5.6. Under the assumptions of this section
$$p\le 11.$$

{\sl Proof.} Assume $p > 11$. From Lemma 5.4 we infer that $\phi:X\to \bbP^1$ has $mp$ 
singular fibres away from the fibre containing $F$.
We use the well-known formula for the Euler characteristic of an elliptic surface
([CD], Proposition 5.1.6):
$$24 = e(X) = c_2(X) = \sum_{b\in \bbP^1}(e(X_b)+\delta(X_b)),\eqno (5.3)$$
where $e(X_b)$ is the Euler (\'etale) characteristic of fibre $X_b$ and 
$\delta(X_b)$ is a certain invariant of wild ramification. 
We have $\delta(X_b) \ge 0$ and $\delta(X_b) = 0$ unless $p = 2,3$ and $X_b$ is
singular of additive type. Under our assumption, $\delta(X_b) = 0$. 
Let $f:V\to \bbP^1$ be the
elliptic fibration of a minimal rational elliptic surface birationally equivalent to
$Y$. 
Let $X_{b_0}$ be the fibre of $\phi$ fixed pointwisely by $g$ and $V_0 = V_{b_0}$ be
the corresponding fibre of $f$. The cover 
$\pi': X\setminus X_{b_0}\to V\setminus V_0$ induced by $\pi:X\to Y$ is unramified and 
moreover splits over each fibre. We have the similar formula for $V$:
$$12 = e(V) = \sum_{b\in \bbP^1}e(V_b).$$ 
We have $e(X_b) = e(\pi(X_b))$ if $X_b \ne X_{b_0}$. 
Thus we get $p(\sum_{b\ne b_0}e(V_b)) \le 24$ and $p > 11$ would imply that 
$\sum_{b\ne b_0}e(V_b) = 1$ and hence $e(V_0) = 11.$ It follows from Kodaira's 
classification of fibres that $V_0$ contains 10 or 11 irreducible components, 
a contradiction to the fact that the Picard number of $V$ is equal to 10.

\med
{\bf Remark 5.7.} In fact, when $p > 3$, one can give a complete
classification of possible degenerate fibres of the elliptic fibration
$\phi:X\to \bbP^1$. We give it in A-D-E notation, recalling the dictionary
for the types of  fibres:
$$I_0 = A_0, I_1 = A_0^*, I_n = A_{n-1} (n > 1), \quad II= A_0^{**},\quad III = A_1^*,
$$
$$IV = A_2^*,\quad I_n^* = D_{n+4},\quad IV^* = E_6,\quad III^* = E_7,\quad II^* =
E_8.$$ 
We have (in each case the type of the fibre $X_{b_0}$ appears first):

\med

$p=11$: 
$$A_0^{**} + 22A_0^* , \quad A_0^{**} + 11A_0^{**};$$

\med

$p=7$:
$$E_8 + 14A_0^*,\quad E_8 + 7A_0^{**},\quad A_1^* + 7A_1^*, \quad A_1^* +
7A_0^* + 7A_0^{**}, \quad A_1^* + 7A_0^* + 7A_1, \quad A_1^* + 21A_0^*.$$

\med

$p=5$: 
$$E_7 + 5A_1^*,\quad E_7 + 5A_0^* + 5A_0^{**},\quad E_7 + 5A_0^* + 5A_1,
\quad E_7 + 15A_0^*,\quad A_2^* + 5A_2^*, \quad A_2^* + 5A_2 + 5A_0^*,$$ 
$$ A_2^* + 5A_1^* + 5A_0^*,\quad A_2^* +10A_0^{**}, \quad A_2^* + 5A_0^{**}+ 5A_1,
\quad A_2^* + 5A_0^{**} + 10A_0^*, \quad A_2^* + 5A_1 + 10A_0^*,\  A_2^* +
20A_0^*.$$

When $p = 3$ and the fibration is elliptic the situation is complicated because of the
presence of wild ramification. However, if the fibration is quasi-elliptic we can have
only the folowing posibilities:
\med
$$A_2^*+9A_2^*,\quad A_2^*+3E_6,\quad E_8+6A_2^*$$
if $p = 3$. In the case $p = 2$ there are many more possibilities. We only indicate
the possible type of the fixed fibre:
$II, III, II^*, I_{2k}^*, 0\le k \le 8.$

\med

To prove the irredundancy of this list (in the case of elliptic fibrations), we provide
an example  for each case.

\med
{\bf Examples 5.8.} In each example below the K3-surface $X$ is given by a
Weierstrass  equation, the automorphism $g$ is given by
$g(x,y,t)=(x,y,t+1)$, and the fixed locus
$X^g$ is the support of the fibre at $t=\infty$.

\med
$p=11 :$

\med
11-(1) $y^2=x^3+x^2+t^{11}-t$

$A_0^{**}$ at $t=\infty$, $11A_0^*$ at $t^{11}-t=0$ and  $11A_0^*$ at $t^{11}-t=-4/27$.

\med
11-(2) $y^2=x^3+t^{11}-t$

$A_0^{**}$ at $t=\infty$ and $11A_0^{**}$ at $t^{11}-t=0$.

\med
$p=7 :$

\med
7-(1) $y^2=x^3+x^2+t^{7}-t$

$E_8$ at $t=\infty$, $7A_0^*$ at $t^{7}-t=0$ and  $7A_0^*$ at $t^{7}-t=-4/27$. 

\med
7-(2) $y^2=x^3+t^{7}-t$

$E_8$ at $t=\infty$ and $7A_0^{**}$ at $t^{7}-t=0$.

\med
7-(3) $y^2=x^3+(t^{7}-t)x$

$A_1^{*}$ at $t=\infty$ and $7A_1^*$ at $t^{7}-t=0$.

\med
7-(4) $y^2=x^3+(t^{7}-t)(x+1)$

$A_1^{*}$ at $t=\infty$, $7A_0^{**}$ at $t^{7}-t=0$ and  $7A_0^*$ at
$t^{7}-t=-4/27$.

\med
7-(5) $y^2=x^3+x^2+(t^{7}-t)x$

$A_1^{*}$ at $t=\infty$, $7A_1$ at $t^{7}-t=0$ and  $7A_0^*$ at $t^{7}-t=1/4$.

\med
7-(6) $y^2=x^3+(t^{7}-t)x+1$

$A_1^{*}$ at $t=\infty$ and $21A_0^*$ at $(t^{7}-t)^3=-27/4$

\med
$p=5 :$

\med
5-(1) $y^2=x^3+(t^{5}-t)x$

$E_7$ at $t=\infty$ and $5A_1^*$ at $t^{5}-t=0$.

\med
5-(2) $y^2=x^3+(t^{5}-t)(x+1)$

$E_7$ at $t=\infty$, $5A_0^{**}$ at $t^{5}-t=0$ and  $5A_0^*$ at
$t^{5}-t=-27/4$.

\med
5-(3) $y^2=x^3+x^2+(t^{5}-t)x$

$E_7$ at $t=\infty$, $5A_1$ at $t^{5}-t=0$ and  $5A_0^*$ at $t^{5}-t=1/4$.

\med
5-(4) $y^2=x^3+(t^{5}-t)x+1$

$E_7$ at $t=\infty$ and $15A_0^*$ at $(t^{5}-t)^3=-27/4$.

\med
5-(5) $y^2=x^3+(t^{5}-t)^2$

$A_2^*$ at $t=\infty$ and $5A_2^*$ at $t^{5}-t=0$.

\med
5-(6) $y^2=x^3+(x+t^{5}-t)^2$

$A_2^*$ at $t=\infty$, $5A_2$ at $t^{5}-t=0$ and  $5A_0^*$ at $t^{5}-t=4/27$.

\med
5-(7) $y^2=x^3+(t^{5}-t)x+(t^{5}-t)^2$

$A_2^*$ at $t=\infty$, $5A_1^*$ at $t^{5}-t=0$ and  $5A_0^*$ at
$t^{5}-t=-4/27$.

\med
5-(8) $y^2=x^3+(t^{5}-t)^2-1$

$A_2^*$ at $t=\infty$ and $10A_0^{**}$ at $(t^{5}-t)^2-1=0$.

\med
5-(9) $y^2=x^3+3x^2+3(t^{5}-t)x+(t^{5}-t)^2$

$A_2^*$ at $t=\infty$, $5A_1$ at $t^{5}-t=0$ and  $5A_0^{**}$ at
$(t^{5}-t)-1=0$.

\med
5-(10) $y^2=x^3+(t^{5}-t)x+(t^{5}-t)(1-t^{5}+t)$

$A_2^*$ at $t=\infty$, $5A_0^{**}$ at $t^{5}-t=0$ and  $10A_0^*$ at
$(t^{5}-t)^2+1=0$.

\med
5-(11) $y^2=x^3+x^2+(t^{5}-t)^2$

$A_2^*$ at $t=\infty$, $5A_1$ at $t^{5}-t=0$ and  $10A_0^*$ at
$(t^{5}-t)^2=-4/27$.

\med
5-(12) $y^2=x^3+x^2+(t^{5}-t)^2-1$

$A_2^*$ at $t=\infty$, $10A_0^*$ at $(t^{5}-t)^2-1=0$ and  $10A_0^*$ at
$(t^{5}-t)^2-1=-4/27$.

\med

 Let us give one application of so far obtained results to
elliptic fibrations with a section on a K3-surface. Assume $f:X\to \bbP^1$ is such a
fibration with a torsion section $s$ of order $p = \char k$. Then the translation
automorphism $g$ is of order $p$ and we can apply our results to it. Let
$s_0$ be the identity section. Let $\calN \to \bbP^1$ be the Neron model of $f$, the 
group scheme obtained by deleting from $X$ the singular locus of every fibre. Let 
$\calN^0$ be its connected component, the largest group subscheme of $\calN$ with
connected fibres. Then the section $s$ is a section of the Neron
model $\calN \to \bbP^1$. If it intersects $\calN^0$ at a fibre of
multiplicative type ($A_0^*$ or $A_n, n > 0$), the fibre must be
pointwisely fixed. If it intersects a fibre of multiplcative type outside
$\calN^0$, then this fibre must be of type $A_{ap-1}$, on which the automorphism $g$
has no fixed points. If $s$ intersects a
fibre of additive type (other singular types) and $p\ne 2,3$, then $s$ is a
section of $\calN^0$ and the locus of singular points of the fibre
belongs to the fixed locus $F$. (Furthermore, $F$ contains central components if $p >
2$ and the fibre is of type $D_n, E_7, E_8$, or if $p > 3$ 
and the fibre is of type $D_n, E_6, E_7, E_8$.)
If $p = 2$, the section  $s$ may intersect a fibre of type $A_1^*,E_7$ or
$D_n$ outside
$\calN_0$. In this case the set of fixed points of $g$ on the fibre is again
non-empty.

If $p = 3$, the section $s$ may intersect a fibre of type $A_2^*,E_6$ outside
$\calN^0$. In this case the set of fixed points of $g$ on this fibre is
equal to an isolated point if the type is $A_2^*$ and to an isolated point on the
central
component if the type is $E_6$.  Since we know that
$F$ is connected and the number of isolated fixed points is at most 2, we
obtain the following information:
\item{(i)} The fibration $f$ has no fibre of type $A_0^*, A_{n}, n\ne ap-1$, if $p > 2$. 
(More generally, if $p > 2$, 
a translation automorphism by a $p$-torsion can not fix a fibre pointwisely, as 
was proved in the proof of Lemma 5.2.) This is still true even if $p=2$. (See Corollary 5.11
below.)
\item{(ii)} If $f$ has a fibre of type $A_0^{**},A_1^*,A_2^*,
$ then it has at most two such fibres, and, if
$p > 3$, then its other possible fibres are of type $A_{ap-1}$;
\item{(iii)} if $f$ has a fibre of type $D_n,E_6,E_7$ or $E_8$, and
$p > 3$, then it has at most one such fibre and all other singular fibres
are of type
$A_{ap-1}$.

\med

\proclaim Corollary 5.9. Assume $p =\char k > 19$. Then an elliptic fibration on a
K3-surface has no $p$-torsion sections.

{\sl Proof.} Use the formula (5.3). In case (ii) from above, we get
$n$ fibres of type $A_{a_ip-1}$  and at most two fibres of type
$A_0^{**},A_1^*,A_2^*$. This gives 
$p(\sum_{i=1}^na_i) \le 20$, and hence $p\le 19$.

In case (iii) from above, we get 
$p(\sum_{i=1}^na_i) \le 24$ and we get $p\le 23$. In fact, if $p = 23$ we
have one fibre of type $A_{22}$ which has two many irreducible components. 
So $p\le 19$ again.

\med
{\bf Remark 5.10.} If $p =2$ and an elliptic fibration has a section, we
have an obvious involution of order 2  which is the
 induced by the inversion automorphism
$x\to -x$ of the general fibre. 
The fixed locus of this involution contains the identity section and any
irreducible fibre of additive type. On each multiplicative fibre of type
$A_{2n}$ it has an isolated fixed point. This looks like a contradiction to
the connectedness of $X^g$. However, it is not. If our results were true that we
sincerely hope so, it implies that in the case of presence of
multiplicative fibres with odd number of irreducible components there is
another one-dimensional part of the fixed locus which intersects the
multiplicative fibres at its fixed point. This is the curve $R$ which is an
inseparable cover of degree 2 of the base which after a base change
defining this cover becomes a section of order 2. Recall that the
Weierstrass equation for an elliptic curve $E$ over a field $K$ of
characteristic 2 looks as follows:
$$y^2+xy = x^3+ax+b \quad\hbox{if $j(E) \ne 0$},$$
$$y^2+cy = x^3+ax+b \quad\hbox{if $j(E) = 0$}.$$
The inversion automorphism is given by $(x,y)\to (x,y+x)$ in the former
case and by $(x,y+c)$ in the latter case. From this we see that the point
$(0,\sqrt{b})$ is the nontrivial 2-torsion point of $E(\bar K)$ if $j(E)\ne
0$ and the group of 2-torsion points is trivial if $j(E)= 0$. In our case
when $K = k(\bbP^1)$, we see that $E(\bar K)_2 =\{0\}$ only when all
fibres are of additive type. Otherwise, either $E(K)_2\ne \{0\}$
or $E(L)_2\ne 0$ for some purely inseparable double cover $L$ of $K$. So we
derive the following:

\med

\proclaim Corollary 5.11.  Let $f:X\to \bbP^1$ be an
elliptic fibration of a K3-surface over a field of characteristic $p = 2$. 
Assume there exists a non-trivial section of order 2. Then all
fibres of multiplicative type contain even number of irreducible
components.

\bigskip

\centerline{{\bf 6. $\kappa(X,F) = 2$.}}

\med
As in the previous case, using the Zariski decomposition, we get that $nF = P'+N'$, 
where $P'$ is a nef divisor with $P'^2 > 0$ and all components of $N'$ are contained
in $P'$. 
This shows that $F$ is equal to the support of a nef divisor $D$ with $D^2 >0$. 
Replacing $D$ by $2D$ we get a linear system without fixed part. 
Also by lemma 3.10, $|2D|$ has no base-points and defines a  morphism of degree $\le
2$ onto 
a surface $X'$ in $\bbP(H^0(X,\calO_X(2D))^*) = \bbP^N$ with at most rational double
points.

\proclaim Theorem 6.1. Let $D$ be the smallest nef divsior with support
on $F$ such that the linear system has no fixed part. Let $N =
{D^2\over 2}+1$ and $d = \dim H^0(X,\calO_X(D-F))$. Then we have the
following inequality
$$p(N-d-1)\le 2N-2.$$

{\sl Proof.} This is a generalization of Oguiso's argument from [Og2],
where this result was proven under the assumption that $F$ is irreducible
(and hence $D = F$). We refer for the details of the proof to Oguiso's
paper and only sketch the main idea. 

Let $\phi_{D}:X\to \bbP(H^0(X,\calO_X(D))^*) = \bbP^N$ be the map given
by the linear system $|D|$. By Lemma 3.10, it is a morphism of degree at
most 2. Thus the image $X'$ of $X$ under $\phi_{|D|}$ is a surface of
degree $D^2 = 2N-2$ or $D^2/2 = N-1$. Since
$D$ is fixed by
$g$ pointwisely, the linear system is $G$-invariant and the induced 
action of $G$ on $X'$ is linear. The exact sequence
$$0\to \calO_X(D-F)\to \calO_X(D)\to \calO_F(D)\to 0$$
shows that the image $F'$ of $F$ spans a linear subspace 
$$H =
\bbP(H^0(X,\calO_X(D-F))^\perp)\subset \bbP(H^0(X,\calO_X(D))^*).$$ It is
easy to see that one can choose a projection from $f:\bbP^N\to H$ given
by a linear sysbsystem $|V|$ where $V =H^0(X,\calO_X(D))^g$ is the linear
subspace of $g$-invariant elements from
$H^0(X,\calO_X(D))$. The image $\bar X'=f(X')$ is an irreducible
non-degenerate surface in
$H$, and hence has degree $\ge \dim H-1 = N-d-1$. On the other hand, the
composition $f\circ \phi_{|D|}:X\to H$ factors through the quotient surface
$Y= X/G$ and hence the degree of the map $f:X'\to \bar X'$ is divisible by
$p$. Thus
$$p(N-d-1) \le p\deg(\bar X') \le \deg X' \le 2N-2.$$
This proves the assertion. 

Under additional assumption that $F$ is irreducible, 
we see that $d = 1$, and, when $p \ne 2$, this easily gives
$(p,N) = (3,2),(3,3),(3,4)$, or $(s,2), s\ge 5$, and $\phi_{|2D|}$ is birational.  
The case $N=2$ (a double plane) requires additional consideration, and Oguiso showed
that $s=5$.

\med
Finally let us reprove, by our means, another result of Oguiso.

\med
\proclaim Theorem 6.2. Assume $F$ is irreducible and $p > 2$. Then $F$ is a
rational curve of arithmetic genus $> 1$.

{\sl Proof.} Consider the branch curve $B\subset Y$. By Corollary 3.3,
$K_Y\sim -(p-1)B$ and, by Corollary 3.6, $H^1(((p-1)B,\calO_{(p-1)B}) = k$.
We know that $F = \pi^*(B)_{red}$. Since $F$ is irreducible, $B$ is also
irreducible. Let $B = a\bar B$ as Weil diviors, where $\bar B$ is
reduced. The restriction of $\pi$ to $F$ defines  a map $\bar \pi:F\to
\bar B$ of degree 1 or purely inseparable of degree $p$. Consider the exact
sequence
$$0\to \calO_Y(-\bar B)/\calO_Y(-(a(p-1)\bar B))\to \calO_{a(p-1)\bar B}\to
\calO_{\bar B}\to 0.\eqno (6.1)$$
Since $\bar B$ is irreducible and reduced, we have $h^0(\calO_{\bar B}) =
1$. Thus if we show that 
$$H^1(Y,\calO_Y(-\bar B)/\calO_Y(-(a(p-1)\bar B)) \ne
0$$
we will be able to infer from (6.1) that $H^1(\bar B,\calO_{\bar B}) = 0$.
This would imply that
$\bar B$ is a smooth rational curve. Since $F$ is (Zariski) homeomeorphic to
$\bar B$, we see that the \'etale Euler characteristics of $F$ and $\bar
B$ are equal. This immediately implies that the normalization of $F$ is a
smooth rational curve. The arithmetic genus of $F$ is of course greater
than 1 since $F^2 > 0$ under the assumption that $\kappa(X,F) = 2$. 

It remains to show that $H^1(Y,\calO_Y(-\bar B)/\calO(-(a(p-1)\bar B)))\ne
0$. We have the obvious exact sequence
$$H^1(Y,\calO_Y(-\bar B))\to H^1(Y,\calO_Y(-\bar B)/\calO_Y(-(a(p-1)\bar
B)))\to $$
$$H^2(Y,\calO_Y(-a(p-1)\bar B)) \to H^2(Y,\calO_Y(-\bar
B))\to 0.$$ Since $\bar B$ is irreducible and $\nu(Y,\bar B) = 2$, we can
apply Theorem (5.1) from [Sa3] to get
$$H^1(Y,\calO_Y(-\bar B)) = 0.$$ By duality 
$$H^2(Y,\calO_Y(-a(p-1)\bar B)) = H^0(Y,\calO_Y(K_Y+(a(p-1)\bar B)) = 
H^0(Y,\calO_Y) = k,$$
$$H^2(Y,\calO_Y(-\bar B)) = H^0(Y,\calO_Y(K_Y+\bar B)) =
H^0(Y,\calO_Y((1-a(p-1))\bar B)) = 0.$$
Here we use that $p > 1$ so that $1-a(p-1) < 0$. This proves the
assertion.

\med
{\bf Remark 6.3.} Assume $p = 2$ and let $f:X\to \bbP^1$ be an elliptic
fibration with a section. As we have explained in Remark 5.10, the fixed
locus $F$ of the inversion involution has the following structure: 

\smallskip\noindent
Case 1 ($j$ is constant): $F$ is the union of the zero section and
the supports of singular fibres (all of them are of additive type). 

\smallskip\noindent
Case 2 ($j$ is not constant)

(a): $F$ is the union of two 2-torsion  sections
and the supports of the fibres of additive type. 

(b): $F$ is the union of the identity section, an irreducible
 curve $C$ which covers the base inseparably of degree 2, and the
supports of the fibres of additive type.

\vfill\eject

\noindent
{\bf References}

\medskip\noindent
[Ar1] M. Artin, Wildly ramified $\bbZ/2$ actions in dimension two, {\it Proc.
A.M.S}, 52 (1975),60-64.

\smallskip\noindent
[Ar2] M. Artin, Coverings of the rational double
points in characteristic
$p$, in {\it Complex Analysis and Algebraic Geometry}, Cambridge Univ.
Press, (1977), 11-22.

\smallskip\noindent
[ABK] A. Assadi, R. Barlow, F. Knop, Abelian group actions on algebraic varieties with
one fixed point, {\it Math. Z.} {\bf 191} (1992), 129-136. 

\smallskip\noindent
[CE] A. Cartan, S. Eilenberg, {\it Homological Algebra}. Princeton Univ.
Press.1956.

\smallskip\noindent
[CD] F. Cossec, I. Dolgachev, {\it Enriques Surfaces I}, Birkh\"auser, 1989.

\smallskip\noindent
[CZ] D. Cox, S. Zucker, Intersection numbers of sections of elliptic surfaces, {\it 
Invent. Math.}, {\bf 53}, (1979), 1-44.

\smallskip\noindent
[G] A. Grothendieck, Sur quelques points d'alg\'ebre homologique, {\it Tohoku Math.
J.} {\bf 9} (1957),
119-221.

\smallskip\noindent
[Ha] R. Hartshorne, {\it Residues and Duality}, Lecture Notes in Math. 20, 
Springer-Verlag, 1966. 

\smallskip\noindent
[La] W. Lang, Extremal  rational elliptic surfaces in characteristic p.
II.: surfaces with three or fewer singular fibres, {\it Arkiv f\"or Matematik},
32 (1994),423-448.

\smallskip\noindent
[Mu] D. Mumford, {\it Geometric Invariant Theory}. Springer-Verlag. 1965

\smallskip\noindent 
[Ni1] V. Nikulin, Finite groups of automorphisms of K\"ahlerian surfaces of
type K3,  {\it Trans. Mosc. Math. Soc.} {\bf 38} (1980), 71-135.

\smallskip\noindent
[Ni2] V. Nikulin, Quotient-groups of groups of automorphisms of hyperbolic 
forms by subgroups generated by $2$-reflections. Algebro-geometric
applications, in `` Current problems in mathematics'', Vol. 18, pp. 3--114,
Akad. Nauk SSSR, Vsesoyuz. Inst. Nauchn. i Tekhn. Informatsii,
Moscow,1981 [English translation: J. Soviet Mathematics, 22
(1983),1401-1476].  

\small\noindent
[Ni3] V. Nikulin, Del Pezzo surfaces with log-terminal singularities. 
Matematicheskii Sbornik, 180 (1989), 226--243.

\smallskip\noindent
[Og1] K. Oguiso, Application of height paring, (in Japanese), {\it Proc.
Minisymposium on algebraic geometry, University of Tokyo}, (1989), 52-66.

\smallskip\noindent
[Og2] K. Oguiso, A note on $\bbZ/p\bbZ$-actions on K3 surfaces in odd characteristic
$p$, {\it Math. Ann.} {\bf 286} (1990), 735-752.

\smallskip\noindent
[Pe1] B. Peskin, On rings of invariants with rational singularities, {\it Proc. Amer. 
Math. Soc.} {\bf 87} no.4 (1983), 621-626.

\smallskip\noindent
[Pe2] B. Peskin, Quotient singularities and wild p-cyclic actions, {\it J. of Algebra},
 81 (1983), 72-99.

\smallskip\noindent
[Pe3] B. Peskin, On the dualizing sheaf of a quotient scheme, {\it Communications in
Algebra}, 12 (1984),1855-1869.

\smallskip\noindent
[Re] M. Reid, Chapters on algebraic surfaces, in ``Complex Algebraic geometry'',
IAS/Park City Math.
Series, vol. 3, 1998.

\smallskip\noindent
[Sa1] F. Sakai, Anticanonical models of rational surfaces, {\it Math. Ann.}
 269 (1984), 389-410.

\smallskip\noindent
[Sa2] F. Sakai, Classification of normal surfaces, in ``Algebraic Geometry.
Bowdoin, Part 1, pp. 451-467. Providence. 1985. 

\smallskip\noindent 
[Sa3] F. Sakai, Weil divisors on normal surfaces, {\it Duke Math. J.} 51 (1984),
877-887.

\smallskip\noindent
[Sh] T. Shioda,  On the Mordell-Weil lattices, {\it Commentarii Mathematici Univ.
Sancti Pauli} 39 no.2(1990), 211-240.

\smallskip\noindent
[Ta] Y. Takeda, Artin-Shreier coverings of algebraic surfaces, {\it J. Math. Soc.
Japan} 
41 (1989), 415-435.

\vglue .3in
\leftline{I. Dolgachev}
\leftline{Department of Mathematics, University of Michigan, Ann Arbor, MI 48109, USA.}
\leftline{e-mail: idolga@math.lsa.umich.edu}

\med
\leftline{J. Keum}
\leftline{Department of Mathematics, Konkuk University, Seoul 143-701, Korea.}
\leftline{e-mail: jhkeum@kkucc.konkuk.ac.kr}

\bye